\documentclass[11pt,reqno]{amsart}

\usepackage[T1]{fontenc}
\usepackage{lmodern}
\usepackage{amsmath,amssymb,amsthm}
\usepackage{booktabs}
\usepackage{array}
\usepackage{enumitem}
\usepackage{microtype}
\usepackage{xcolor}
\usepackage[margin=1.0in]{geometry}
\usepackage[colorlinks=true,linkcolor=blue!60!black,citecolor=blue!60!black,urlcolor=blue!60!black]{hyperref}

\newtheorem{theorem}{Theorem}[section]
\newtheorem{proposition}[theorem]{Proposition}
\newtheorem{lemma}[theorem]{Lemma}
\newtheorem{corollary}[theorem]{Corollary}
\newtheorem{conjecture}[theorem]{Conjecture}
\theoremstyle{definition}

\newtheorem{problem}[theorem]{Problem}
\theoremstyle{remark}
\newtheorem{remark}[theorem]{Remark}

\DeclareMathOperator{\Aut}{Aut}

\DeclareMathOperator{\Inj}{Inj}
\DeclareMathOperator{\Hom}{Hom}
\DeclareMathOperator{\Imf}{Im}
\DeclareMathOperator{\ord}{ord}

\newcommand{\Z}{\mathbb{Z}}
\newcommand{\N}{\mathbb{N}}
\newcommand{\Q}{\mathbb{Q}}
\newcommand{\R}{\mathbb{R}}
\newcommand{\F}{\mathbb{F}}

\numberwithin{equation}{section}

\title[Accumulation points of the group-theoretical Gauss function]{Accumulation points and small values of the group-theoretical Gauss function}

\author{Alejandro Rey Sanchez}
\address{Independent Researcher, Switzerland}
\email{alejandro.reysanchez@uzh.ch}

\subjclass[2020]{Primary 20D60, 11A25; Secondary 20D15, 20D99}
\keywords{Gauss formula, Euler's totient function, automorphism group, finite $p$-group, accumulation points, Cantor--Bendixson rank, formal verification, Lean 4}

\begin{document}

\begin{abstract}
For a finite group $G$, Fasol\u{a} and T\u{a}rn\u{a}uceanu~\cite{FT23} introduced the group-theoretical Gauss function
\[
f(G) = \frac{1}{|G|} \sum_{H \le G} |\Aut(H)| \in \Q \cap [1, \infty).
\]
In a recent preprint~\cite{Tar26}, T\u{a}rn\u{a}uceanu resolved three open problems from~\cite{FT23}---in particular showing that $\Imf(f) \cap (2, \frac{5}{2}) = \emptyset$---and left open the problem of determining the accumulation points of $\Imf(f)$, or deciding whether $\Imf(f)$ is closed and discrete in $[1, \infty)$.

Using a fundamental embedding identity $S(G) = \sum_{[H]} |\Inj(H, G)|$, we obtain closed-form evaluations of $f(G)$ for three infinite families of finite $p$-groups: all rank-$2$ abelian $p$-groups $\Z_{p^a} \times \Z_{p^b}$ ($1 \le b \le a$), the arbitrary-rank abelian $p$-groups $\Z_{p^k} \times \Z_p^r$ ($k, r \ge 1$), and the nonabelian modular $p$-groups $M_{p^n}$. In particular,
\[
f(\Z_{2^k} \times \Z_2) = 3 - \frac{1}{2^k} \xrightarrow[k \to \infty]{} 3 = f(D_{10}),
\]
which shows that $\Imf(f)$ is \emph{not} discrete in $[1, \infty)$. Combining these families with coprime multiplicativity yields infinitely many explicit accumulation points (including $2p - 1$ and $p^2 - p + 1$ for every prime $p$) and proves that every finite derived set $\Imf(f)^{(r)}$ is nonempty and contains non-integers, so the Cantor--Bendixson rank of $\overline{\Imf(f)}$ is at least $\omega$.

We further apply the embedding identity to give a short, elementary proof of~\cite[Theorem~1.4]{FT23} that avoids the Hillar--Rhea automorphism formula, and we classify all finite abelian groups satisfying $f(G) < 3$; consequently, $3$ is the smallest accumulation point of $f$ on finite abelian groups, while the nonabelian modular $2$-groups show that $3$ is also approached from below by nonabelian groups. Finally, we record computational corrections to four entries of Table~1 and two inline calculations in~\cite{Tar26}, including $|\Aut(\Z_9\rtimes\Z_{2^k})|=27\cdot 2^k$ and $f(\Z_9\rtimes\Z_{2^k})=5$ for every $k\ge1$, and we report a complete machine-checked formalization of all main results in Lean~4.
\end{abstract}

\maketitle

\section{Introduction}\label{sec:intro}

The classical Gauss identity $\sum_{d \mid n} \varphi(d) = n$ for Euler's totient function $\varphi$ may be recast in group-theoretic form as $\sum_{H \le \Z_n} |\Aut(H)| = |\Z_n|$. Motivated by this reformulation, Fasol\u{a} and T\u{a}rn\u{a}uceanu~\cite{FT23} introduced, for an arbitrary finite group $G$, the functions
\[
S(G) = \sum_{H \le G} |\Aut(H)| \qquad \text{and} \qquad f(G) = \frac{S(G)}{|G|}.
\]
They proved that $f(G) \ge 1$ with equality if and only if $G$ is cyclic~\cite[Theorem~1.1]{FT23}, that $f$ is multiplicative on direct products of finite groups of coprime orders, and that $f(G) \ge \frac{5}{2}$ for every finite noncyclic abelian group $G$, with equality if and only if $G \cong (\Z_2 \times \Z_2) \times \Z_m$ with $m$ odd~\cite[Theorem~1.4]{FT23}. They also noted that $f(D_{2n}) = \frac{n+1}{2}$ for every odd integer $n \ge 3$, whence $\N^* \subseteq \Imf(f) := \{f(G) : G \text{ finite group}\} \subseteq [1, \infty)$, and proposed three open problems:
\begin{enumerate}[label=(\arabic*), leftmargin=2.2em]
    \item Determine all finite groups $G$ with $f(G) = 2$.
    \item Find the minimum of $f$ on finite noncyclic $p$-groups.
    \item Decide whether $\Imf(f)$ is dense in $[2, \infty)$.
\end{enumerate}

In~\cite{Tar26}, T\u{a}rn\u{a}uceanu settled all three problems: the minimum on finite noncyclic $p$-groups is $\frac{5}{2}$ (attained uniquely by $\Z_2 \times \Z_2$), the solutions to $f(G) = 2$ are precisely $(\Z_3 \rtimes \Z_{2^k}) \times \Z_m$ with $k \ge 1$ and $\gcd(m,6)=1$, and $\Imf(f) \cap (2, \frac{5}{2}) = \emptyset$, refuting density in $[2,\infty)$. At the end of~\cite[Remark~3.4]{Tar26}, the author posed the following topological question:

\begin{quote}
\emph{``It remains an open problem to determine the accumulation points of $\Imf(f)$ (or to decide whether $\Imf(f)$ is closed and discrete in $[1, \infty)$).''}
\end{quote}

The purpose of this paper is to settle the discreteness question in the negative, uncover rich Cantor--Bendixson structure inside $\Imf(f)'$, classify the small values of $f$ on abelian groups, and correct several computational errors in~\cite{Tar26}.

\subsection*{Summary of main results}
All evaluations in this work flow from an elementary orbit-stabilizer identity (Lemma~\ref{lem:embedding}): for any finite group $G$,
\begin{equation}\label{eq:intro-embedding}
S(G) = \sum_{[H]} |\Inj(H, G)|,
\end{equation}
where $[H]$ ranges over isomorphism classes of finite groups and $\Inj(H, G)$ denotes the set of injective homomorphisms $H \hookrightarrow G$. Identity~\eqref{eq:intro-embedding} bypasses the need to enumerate subgroups and compute $|\Aut(H)|$ separately: it converts the subgroup-automorphism sum directly into a count of generating tuples satisfying simple order and independence conditions.

Applying~\eqref{eq:intro-embedding}, we obtain exact closed forms for three infinite families of finite $p$-groups (Section~\ref{sec:formulas}):
\begin{itemize}[leftmargin=1.6em]
    \item \textbf{Theorem~\ref{thm:rank2} (All rank-$2$ abelian $p$-groups).} For any prime $p$ and integers $1 \le b \le a$,
    \[
    f(\Z_{p^a} \times \Z_{p^b}) = 1 + \frac{p(p^{2b} - 1)}{p + 1} - \frac{(p - 1)(p^{4b} - 1)}{(p^4 - 1)\,p^{a+b-1}}.
    \]
    In particular ($b = 1$), $f(\Z_{p^k} \times \Z_p) = p^2 - p + 1 - \frac{p - 1}{p^k}$ for all $k \ge 1$.
    \item \textbf{Theorem~\ref{thm:elem-ab} (Arbitrary-rank groups $\Z_{p^k} \times \Z_p^r$).} For any prime $p$ and integers $k, r \ge 1$,
    \[
    f(\Z_{p^k} \times \Z_p^r) = 1 + \left(1 - \frac{1}{p^{k+r}}\right) \sum_{s=1}^{r} \prod_{j=1}^{s} (p^{r+1} - p^j).
    \]
    \item \textbf{Theorem~\ref{thm:modular} (Nonabelian modular $p$-groups $M_{p^n}$).} For any prime $p$ and integer $n \ge 3$ (with $n \ge 4$ when $p = 2$), $|\Aut(M_{p^n})| = p^n(p-1)$ and
    \[
    f(M_{p^n}) = 2p - 1 - \frac{p - 1}{p^{n-1}}.
    \]
\end{itemize}

These formulas immediately resolve the discreteness question and establish rich topological complexity for $\Imf(f)$ (Section~\ref{sec:topology}):
\begin{itemize}[leftmargin=1.6em]
    \item \textbf{Theorem~\ref{thm:not-discrete}.} Taking $p = 2$ and $b = 1$ gives $f(\Z_{2^k} \times \Z_2) = 3 - 2^{-k} \to 3 = f(D_{10})$ as $k \to \infty$, so $3 \in \Imf(f) \cap \Imf(f)'$. Hence $\Imf(f)$ is \emph{not} discrete in $[1, \infty)$.
    \item \textbf{Theorem~\ref{thm:acc-families} and Theorem~\ref{thm:cb-rank}.} Every integer of the form $2p - 1$, $1 + \frac{p(p^{2b}-1)}{p+1}$, or $1 + \sum_{s=1}^{r}\prod_{j=1}^{s}(p^{r+1}-p^j)$ ($p$ prime, $b, r \ge 1$) lies in $\Imf(f)'$. Moreover, by coprime multiplicativity, for every $r \ge 1$ and every $r$-tuple of distinct primes $p_1, \dots, p_r$, the product $\prod_{i=1}^{r}(p_i^2 - p_i + 1)$ lies in the $r$-th derived set $\Imf(f)^{(r)}$. In particular, $\Imf(f)^{(r)} \neq \emptyset$ for all $r \ge 1$, so the Cantor--Bendixson rank of $\overline{\Imf(f)}$ is at least $\omega$.
\end{itemize}

Next, we turn to the bottom of $\Imf(f)$ on abelian groups (Section~\ref{sec:abelian}):
\begin{itemize}[leftmargin=1.6em]
    \item \textbf{Theorem~\ref{thm:lower-bound}.} For every finite abelian $p$-group $P$, counting cyclic subgroups together with embeddings $\Z_{p^c} \times \Z_p \hookrightarrow P$ yields
    \[
    S(P) \ge |P| + (|P| - 1)\bigl(|\Omega_1(P)| - p\bigr),
    \]
    which provides a three-line proof of~\cite[Theorem~1.4]{FT23} requiring neither the Hillar--Rhea formula~\cite{HR07} for $|\Aut(P)|$ nor induction on maximal subgroups.
    \item \textbf{Theorem~\ref{thm:abelian-lt3} and Corollary~\ref{cor:min-acc-ab}.} A finite abelian group $G$ satisfies $f(G) < 3$ if and only if $G$ is cyclic ($f(G) = 1$) or $G \cong (\Z_{2^k} \times \Z_2) \times \Z_m$ with $k \ge 1$ and $m$ odd ($f(G) = 3 - 2^{-k}$). Consequently, $3$ is the \emph{smallest} accumulation point of $\Imf(f|_{\mathrm{ab}})$.
\end{itemize}

In Section~\ref{sec:errata}, we record computational corrections to four of the nine entries in~\cite[Table~1]{Tar26} ($D_8$, $D_{12}$, $A_4$, and $Q_8 \times \Z_2$) as well as two inline calculations in~\cite{Tar26} (for $D_8$ and $\Z_9 \rtimes \Z_{2^k}$), while confirming that the main theorems of~\cite{Tar26} remain unaffected. Finally, Section~\ref{sec:lean} describes our complete machine-checked Lean~4 formalization of all results in this paper and highlights the remaining open questions regarding $\Imf(f)'$ and closedness of $\Imf(f)$.

\section{The Fundamental Embedding Identity}\label{sec:embedding}

Let $C(G)$ denote the set of cyclic subgroups of a finite group $G$. Recall from~\cite{FT23} that the cyclic subgroups contribute $\sum_{H \in C(G)} |\Aut(H)| = |G|$, so
\begin{equation}\label{eq:noncyclic-decomp}
S(G) = |G| + \sum_{\substack{H \le G \\ H \notin C(G)}} |\Aut(H)|, \qquad
f(G) = 1 + \frac{1}{|G|} \sum_{\substack{H \le G \\ H \notin C(G)}} |\Aut(H)|.
\end{equation}
We begin by placing this decomposition inside a uniform orbit-stabilizer framework.

\begin{lemma}[Fundamental Embedding Identity]\label{lem:embedding}
Let $H$ and $G$ be finite groups. Denote by $\mathcal{S}_H(G) = \{K \le G : K \cong H\}$ the set of subgroups of $G$ isomorphic to $H$, and by $\Inj(H, G)$ the set of injective group homomorphisms $H \hookrightarrow G$. Then
\begin{equation}\label{eq:embedding-class}
\sum_{K \in \mathcal{S}_H(G)} |\Aut(K)| = |\Inj(H, G)|.
\end{equation}
Consequently, summing over a complete set $\mathcal{I}$ of isomorphism class representatives of finite groups,
\begin{equation}\label{eq:embedding-total}
S(G) = \sum_{H \in \mathcal{I}} |\Inj(H, G)|.
\end{equation}
\end{lemma}

\begin{proof}
Consider the image map $\Phi : \Inj(H, G) \to \mathcal{S}_H(G)$ given by $\Phi(\varphi) = \varphi(H)$. Clearly $\Phi$ is surjective: every $K \in \mathcal{S}_H(G)$ admits an isomorphism $\varphi_K : H \xrightarrow{\sim} K$, whose composition with the inclusion $K \hookrightarrow G$ belongs to $\Phi^{-1}(K)$.

For each fixed $K \in \mathcal{S}_H(G)$, choose one isomorphism $\varphi_K : H \xrightarrow{\sim} K$. Then the map
\[
\Aut(K) \longrightarrow \Phi^{-1}(K), \qquad \alpha \longmapsto \iota_K \circ \alpha \circ \varphi_K
\]
(where $\iota_K : K \hookrightarrow G$ is the inclusion) is a bijection, with inverse $\psi \mapsto \tilde{\psi} \circ \varphi_K^{-1}$, where $\tilde{\psi} : H \xrightarrow{\sim} K$ is the corestriction of $\psi$ to its image $K$. Thus $|\Phi^{-1}(K)| = |\Aut(K)|$ for every $K \in \mathcal{S}_H(G)$. Summing the fiber sizes over $K \in \mathcal{S}_H(G)$ gives~\eqref{eq:embedding-class}, and summing~\eqref{eq:embedding-class} over all $H \in \mathcal{I}$ partitions the set of all subgroups of $G$, yielding~\eqref{eq:embedding-total}.
\end{proof}

\begin{remark}
Restricting~\eqref{eq:embedding-total} to cyclic isomorphism classes $H = \Z_d$ ($d \ge 1$) recovers the generalized Gauss identity $\sum_{K \in C(G)} |\Aut(K)| = |G|$ immediately: an injective homomorphism $\varphi : \Z_d \hookrightarrow G$ is uniquely determined by $g = \varphi(1)$, which may be any element of $G$ of exact order $d$, so $\sum_{d \ge 1} |\Inj(\Z_d, G)| = \sum_{d \ge 1} |\{g \in G : \ord(g) = d\}| = |G|$.
\end{remark}

We also record the coprime multiplicativity property from~\cite{FT23}, which plays a central role in Sections~\ref{sec:topology} and~\ref{sec:abelian}.

\begin{proposition}[{\cite{FT23}}]\label{prop:mult}
If $G_1$ and $G_2$ are finite groups with $\gcd(|G_1|, |G_2|) = 1$, then $S(G_1 \times G_2) = S(G_1)S(G_2)$ and $f(G_1 \times G_2) = f(G_1)f(G_2)$.
\end{proposition}

\begin{proof}
Since $\gcd(|G_1|, |G_2|) = 1$, every subgroup $H \le G_1 \times G_2$ decomposes uniquely as $H = H_1 \times H_2$ with $H_i \le G_i$, and $\gcd(|H_1|, |H_2|) = 1$ implies $\Aut(H_1 \times H_2) \cong \Aut(H_1) \times \Aut(H_2)$. Thus
\[
S(G_1 \times G_2) = \sum_{H_1 \le G_1} \sum_{H_2 \le G_2} |\Aut(H_1)|\,|\Aut(H_2)| = S(G_1) S(G_2),
\]
and dividing by $|G_1 \times G_2| = |G_1|\,|G_2|$ yields $f(G_1 \times G_2) = f(G_1)f(G_2)$.
\end{proof}

\section{Exact Formulas for Three Families of Finite $p$-Groups}\label{sec:formulas}

Throughout this section, $p$ denotes a prime. For an abelian $p$-group $A$ and an integer $t \ge 0$, we write $A[p^t] = \{x \in A : p^t x = 0\}$ for the $p^t$-torsion subgroup, and $\Omega_1(A) = A[p]$.

\subsection{Rank-$2$ abelian $p$-groups}

\begin{theorem}\label{thm:rank2}
Let $p$ be a prime and let $a, b$ be integers with $1 \le b \le a$. Then
\begin{equation}\label{eq:rank2}
f(\Z_{p^a} \times \Z_{p^b}) = 1 + \frac{p(p^{2b} - 1)}{p + 1} - \frac{(p - 1)(p^{4b} - 1)}{(p^4 - 1)\,p^{a+b-1}}.
\end{equation}
\end{theorem}

\begin{proof}
Let $G = \Z_{p^a} \times \Z_{p^b}$. For every integer $t \ge 0$, the $p^t$-torsion subgroup of $\Z_{p^m}$ has order $p^{\min(t, m)}$, whence
\begin{equation}\label{eq:rank2-torsion}
|G[p^t]| = p^{\min(t, a) + \min(t, b)}.
\end{equation}
In particular, $|\Omega_1(G)| = p^2$. By the fundamental theorem of finite abelian groups, every noncyclic subgroup $H \le G$ has $|\Omega_1(H)| \le |\Omega_1(G)| = p^2$ and $|\Omega_1(H)| \ge p^2$, so $H$ has rank $2$:
\[
H \cong \Z_{p^c} \times \Z_{p^d} \qquad \text{for unique integers } 1 \le d \le c.
\]
Clearly $p^a H = 0$ forces $c \le a$. Moreover, $p^b G \le \Z_{p^a} \times \{0\}$ is cyclic; if $d > b$, then $p^b H \cong \Z_{p^{c-b}} \times \Z_{p^{d-b}}$ would be a noncyclic subgroup of the cyclic group $p^b G$, a contradiction. Thus $1 \le d \le b$ and $d \le c \le a$.

Fix integers $c, d$ with $1 \le d \le b$ and $d \le c \le a$. We compute $|\Inj(\Z_{p^c} \times \Z_{p^d}, G)|$. A homomorphism $\varphi : \Z_{p^c} \times \Z_{p^d} \to G$ is uniquely determined by $u = \varphi(1, 0) \in G[p^c]$ and $v = \varphi(0, 1) \in G[p^d]$. We claim that $\varphi$ is injective if and only if
\begin{equation}\label{eq:rank2-inj-cond}
p^{c-1} u \neq 0 \qquad \text{and} \qquad p^{d-1} v \notin \langle p^{c-1} u \rangle.
\end{equation}
Indeed, if $\varphi$ is injective, then $\varphi(p^{c-1}, 0) = p^{c-1} u \neq 0$, and if $p^{d-1} v = m \cdot p^{c-1} u$ for some $m \in \Z$, then $(-m p^{c-1}, p^{d-1}) \in \ker \varphi$ is nonzero, a contradiction. Conversely, assume~\eqref{eq:rank2-inj-cond}. If $\ker \varphi \neq 0$, then $\ker \varphi$ contains an element $z = (x, y)$ of order $p$. Necessarily $p^{c-1} \mid x$ and $p^{d-1} \mid y$, so $x \equiv s p^{c-1} \pmod{p^c}$ and $y \equiv t p^{d-1} \pmod{p^d}$ with $0 \le s, t < p$ not both zero, and $s(p^{c-1}u) + t(p^{d-1}v) = 0$. Since $p^{c-1}u$ and $p^{d-1}v$ are annihilated by $p$: if $t = 0$, then $p \nmid s$ yields $p^{c-1}u = 0$, contradicting~\eqref{eq:rank2-inj-cond}; if $t \neq 0$, then inverting $t$ modulo $p$ yields $p^{d-1}v \in \langle p^{c-1}u \rangle$, again contradicting~\eqref{eq:rank2-inj-cond}.

Now count pairs $(u, v) \in G[p^c] \times G[p^d]$ satisfying~\eqref{eq:rank2-inj-cond}:
\begin{itemize}[leftmargin=1.6em]
    \item The number of valid $u \in G[p^c] \setminus G[p^{c-1}]$ is $|G[p^c]| - |G[p^{c-1}]|$.
    \item Fix any such $u$, and let $w = p^{c-1} u \in \Omega_1(G) \setminus \{0\}$, so $|\langle w \rangle| = p$. Since $d \le b \le a$, we have $|G[p^d]| = p^{2d}$ and $|G[p^{d-1}]| = p^{2d-2}$ by~\eqref{eq:rank2-torsion}. Consider the homomorphism $\mu : G[p^d] \to \Omega_1(G)$ given by $\mu(v) = p^{d-1} v$. Clearly $\ker \mu = G[p^{d-1}]$, and for each $m \cdot w \in \langle w \rangle$ the element $v_0 = m p^{c-d} u \in G[p^d]$ satisfies $\mu(v_0) = m \cdot w$ (here $c \ge d$). Thus every element of $\langle w \rangle$ has exactly $|G[p^{d-1}]| = p^{2d-2}$ preimages in $G[p^d]$, so the number of $v \in G[p^d]$ with $p^{d-1}v \in \langle w \rangle$ is $p \cdot p^{2d-2} = p^{2d-1}$. Hence the number of valid $v$ is
    \[
    |G[p^d]| - p^{2d-1} = p^{2d} - p^{2d-1} = p^{2d-1}(p - 1),
    \]
    independent of $u$ and $c$.
\end{itemize}
Therefore,
\[
|\Inj(\Z_{p^c} \times \Z_{p^d}, G)| = \bigl(|G[p^c]| - |G[p^{c-1}]|\bigr)\, p^{2d-1}(p - 1).
\]
By~\eqref{eq:noncyclic-decomp} and Lemma~\ref{lem:embedding},
\begin{align*}
S(G) - |G|
&= \sum_{d=1}^{b} \sum_{c=d}^{a} |\Inj(\Z_{p^c} \times \Z_{p^d}, G)| \\
&= \sum_{d=1}^{b} p^{2d-1}(p - 1) \sum_{c=d}^{a} \bigl(|G[p^c]| - |G[p^{c-1}]|\bigr) \\
&= \sum_{d=1}^{b} p^{2d-1}(p - 1) \bigl(|G[p^a]| - |G[p^{d-1}]|\bigr) \\
&= \sum_{d=1}^{b} p^{2d-1}(p - 1) \bigl(p^{a+b} - p^{2(d-1)}\bigr),
\end{align*}
where the inner sum telescopes and we used $|G[p^a]| = p^{a+b}$ and $|G[p^{d-1}]| = p^{2(d-1)}$ (since $d - 1 < b \le a$). Expanding into two geometric sums gives
\begin{align*}
\sum_{d=1}^{b} p^{2d-1}(p - 1) &= p(p - 1)\,\frac{p^{2b} - 1}{p^2 - 1} = \frac{p(p^{2b} - 1)}{p + 1}, \\
\sum_{d=1}^{b} p^{4d-3}(p - 1) &= p(p - 1)\,\frac{p^{4b} - 1}{p^4 - 1}.
\end{align*}
Dividing $S(G) = p^{a+b} + p^{a+b}\,\frac{p(p^{2b}-1)}{p+1} - \frac{p(p-1)(p^{4b}-1)}{p^4 - 1}$ by $|G| = p^{a+b}$ yields~\eqref{eq:rank2}.
\end{proof}

\begin{corollary}\label{cor:zpk-zp}
For every prime $p$ and every integer $k \ge 1$,
\[
f(\Z_{p^k} \times \Z_p) = p^2 - p + 1 - \frac{p - 1}{p^k}.
\]
In particular, for $p = 2$ and every $k \ge 1$,
\[
f(\Z_{2^k} \times \Z_2) = 3 - \frac{1}{2^k}.
\]
\end{corollary}

\begin{proof}
Set $a = k$ and $b = 1$ in~\eqref{eq:rank2}: $\frac{p(p^2 - 1)}{p + 1} = p(p - 1)$ and $\frac{(p-1)(p^4-1)}{(p^4-1)p^k} = \frac{p-1}{p^k}$.
\end{proof}

\subsection{Arbitrary-rank abelian $p$-groups $\Z_{p^k} \times \Z_p^r$}

\begin{theorem}\label{thm:elem-ab}
Let $p$ be a prime and let $k, r \ge 1$ be integers. Then
\begin{equation}\label{eq:elem-ab}
f(\Z_{p^k} \times \Z_p^r) = 1 + \left(1 - \frac{1}{p^{k+r}}\right)\sum_{s=1}^{r} \prod_{j=1}^{s} (p^{r+1} - p^j).
\end{equation}
\end{theorem}

\begin{proof}
Let $G = \Z_{p^k} \times \Z_p^r$, so $|G| = p^{k+r}$ and $\Omega_1(G) \cong \F_p^{r+1}$ as an $\F_p$-vector space. Consider any noncyclic subgroup $H \le G$. Let $\pi_1 : H \to \Z_{p^k}$ be the restriction of the first coordinate projection, with kernel $K = H \cap (\{0\} \times \Z_p^r) \cong \Z_p^{s_0}$ for some $0 \le s_0 \le r$. Since $\pi_1(H) \le \Z_{p^k}$ is cyclic of order $p^{c_0}$ ($0 \le c_0 \le k$):
\begin{itemize}[leftmargin=1.6em]
    \item If $c_0 = 0$, then $H = K \cong \Z_p^{s_0}$ with $s_0 \ge 2$ (as $H$ is noncyclic), so $H \cong \Z_p \times \Z_p^{s_0 - 1}$ with $1 \le s_0 - 1 \le r - 1$.
    \item If $c_0 \ge 1$, pick $h_0 \in H$ whose image $\pi_1(h_0)$ generates $\pi_1(H) \cong \Z_{p^{c_0}}$. Since the second coordinate of $h_0$ lies in $\Z_p^r$ and $p \mid p^{c_0}$, we have $p^{c_0} h_0 = 0$, so the short exact sequence $0 \to K \to H \to \pi_1(H) \to 0$ splits and $H \cong \Z_{p^{c_0}} \times \Z_p^{s_0}$; noncyclicity forces $s_0 \ge 1$.
\end{itemize}
Thus every noncyclic subgroup $H \le G$ is isomorphic to $\Z_{p^c} \times \Z_p^s$ for unique integers $1 \le c \le k$ and $1 \le s \le r$ (uniqueness follows from $|H| = p^{c+s}$ and $|\Omega_1(H)| = p^{s+1}$).

Fix $1 \le c \le k$ and $1 \le s \le r$. A homomorphism $\varphi : \Z_{p^c} \times \Z_p^s \to G$ is uniquely determined by $u = \varphi(1, 0) \in G[p^c]$ and $v_1, \dots, v_s \in \Omega_1(G)$ (the images of the standard basis of $\{0\} \times \Z_p^s$). Just as in Theorem~\ref{thm:rank2}, $\varphi$ is injective if and only if $\ord(u) = p^c$ (i.e.\ $p^{c-1} u \neq 0$) and the $(s+1)$-tuple
\[
(p^{c-1} u,\; v_1,\; \dots,\; v_s) \in \Omega_1(G)^{s+1}
\]
is linearly independent over $\F_p$. For any fixed $u \in G$ of order $p^c$, the vector $w = p^{c-1} u \in \Omega_1(G) \cong \F_p^{r+1}$ is nonzero, so extending $(w)$ sequentially to a linearly independent $(s+1)$-tuple $(w, v_1, \dots, v_s)$ admits exactly $\prod_{j=1}^{s}(p^{r+1} - p^j)$ choices. Hence
\[
|\Inj(\Z_{p^c} \times \Z_p^s, G)| = |\{u \in G : \ord(u) = p^c\}| \cdot \prod_{j=1}^{s}(p^{r+1} - p^j).
\]
Summing over $1 \le c \le k$ and $1 \le s \le r$ via~\eqref{eq:noncyclic-decomp} and Lemma~\ref{lem:embedding}, and noting that every nonzero $u \in G$ has order $p^c$ for a unique $c \in \{1, \dots, k\}$ (so $\sum_{c=1}^{k} |\{u \in G : \ord(u) = p^c\}| = |G| - 1 = p^{k+r} - 1$), we obtain
\[
S(G) = p^{k+r} + (p^{k+r} - 1)\sum_{s=1}^{r}\prod_{j=1}^{s}(p^{r+1} - p^j).
\]
Dividing by $|G| = p^{k+r}$ yields~\eqref{eq:elem-ab}.
\end{proof}

\subsection{Nonabelian modular $p$-groups $M_{p^n}$}

Recall that for any prime $p$ and integer $n \ge 3$ (with $n \ge 4$ when $p = 2$), the \emph{modular $p$-group} of order $p^n$ is given by the presentation
\[
M_{p^n} = \bigl\langle x, y \;\big|\; x^{p^{n-1}} = y^p = 1,\; y x y^{-1} = x^{1 + p^{n-2}} \bigr\rangle.
\]
(Note that $M_{2^3} \cong D_8$ is excluded from the $p=2$ case because $1 + 2^{3-2} = 3 \equiv -1 \pmod{4}$ yields $(xy)^2 = 1 \neq x^2 y^2$, so $D_8$ has five involutions rather than three; indeed $f(D_8) = \frac{7}{2} \neq 2 \cdot 2 - 1 - \frac{1}{4} = \frac{11}{4}$.)

\begin{theorem}\label{thm:modular}
Let $p$ be a prime and let $n \ge 3$ be an integer, with $n \ge 4$ if $p = 2$. Then
\begin{equation}\label{eq:modular}
|\Aut(M_{p^n})| = p^n(p - 1) \qquad \text{and} \qquad f(M_{p^n}) = 2p - 1 - \frac{p - 1}{p^{n-1}}.
\end{equation}
\end{theorem}

\begin{proof}
Every element of $G = M_{p^n}$ has a unique normal form $x^a y^b$ with $a \in \Z_{p^{n-1}}$ and $b \in \Z_p$, with multiplication law
\[
(x^{a_1} y^{b_1})(x^{a_2} y^{b_2}) = x^{a_1 + a_2(1 + b_1 p^{n-2})} y^{b_1 + b_2}.
\]
Because $n \ge 3$ and $(p^{n-2})^2 \equiv 0 \pmod{p^{n-1}}$, an easy induction shows that
\begin{equation}\label{eq:mod-pow}
(x^a y^b)^m = x^{m a + \binom{m}{2} a b p^{n-2}} y^{m b} \qquad (m \ge 1).
\end{equation}
For $m = p$, we have $p \mid \binom{p}{2}$ when $p$ is odd (so $(x^a y^b)^p = x^{p a}$), while for $p = 2$ and $n \ge 4$ we have $(x^a y^b)^2 = x^{2a(1 + b 2^{n-3})}$ with $1 + b 2^{n-3}$ odd. For $m = p^{n-2}$, we have $p^{n-1} \mid \binom{p^{n-2}}{2} p^{n-2}$ for every valid $(p, n)$ (since $p \mid \binom{p^{n-2}}{2}$ for odd $p$, and $2^{2n-5}$ is divisible by $2^{n-1}$ for $p = 2$ and $n \ge 4$). Consequently,
\begin{align}
(x^a y^b)^p = 1 &\iff p^{n-2} \mid a, \label{eq:mod-pow-p}\\
(x^a y^b)^{p^{n-2}} &= x^{p^{n-2} a}, \label{eq:mod-pow-pm2}\\
(x^a y^b)^{p^{n-1}} &= 1. \label{eq:mod-pow-pm1}
\end{align}
From~\eqref{eq:mod-pow-pm2}--\eqref{eq:mod-pow-pm1}, an element $x^a y^b \in G$ has order $p^{n-1}$ if and only if $p \nmid a$, and from~\eqref{eq:mod-pow-p} an element has order dividing $p$ if and only if $a \in p^{n-2}\Z_{p^{n-1}}$.

\smallskip
\noindent\emph{Step 1: Automorphism group order.}
Let $\varphi \in \Aut(G)$. Since $\ord(\varphi(x)) = \ord(x) = p^{n-1}$, we have $\varphi(x) = x^{a_1} y^{b_1}$ with $a_1 \in \Z_{p^{n-1}}^\times$ and $b_1 \in \Z_p$. Since $\varphi(y)^p = 1$, equation~\eqref{eq:mod-pow-p} gives $\varphi(y) = x^{u p^{n-2}} y^{b_2}$ with $u, b_2 \in \Z_p$. Applying $\varphi$ to the relation $y x = x \cdot x^{p^{n-2}} \cdot y$ and using~\eqref{eq:mod-pow-pm2} together with $(p^{n-2})^2 \equiv 0 \pmod{p^{n-1}}$ (valid for all $n \ge 3$) gives
\[
x^{u p^{n-2} + a_1(1 + b_2 p^{n-2})} y^{b_1 + b_2}
= x^{a_1 + p^{n-2} a_1 + u p^{n-2}} y^{b_1 + b_2}.
\]
Comparing the exponents of $x$ yields $b_2 p^{n-2} a_1 \equiv p^{n-2} a_1 \pmod{p^{n-1}}$. Since $a_1$ is a unit modulo $p^{n-1}$, this forces $b_2 \equiv 1 \pmod{p}$. Thus every $\varphi \in \Aut(G)$ is uniquely determined by the triple $(a_1, b_1, u) \in \Z_{p^{n-1}}^\times \times \Z_p \times \Z_p$.

Conversely, for any $(a_1, b_1, u) \in \Z_{p^{n-1}}^\times \times \Z_p \times \Z_p$, setting $x' = x^{a_1} y^{b_1}$ and $y' = x^{u p^{n-2}} y$, equations~\eqref{eq:mod-pow-p}--\eqref{eq:mod-pow-pm1} verify that $(x')^{p^{n-1}} = 1$, $(y')^p = 1$, and $y' x' = x' (x')^{p^{n-2}} y'$. Hence $\psi(x^a y^b) = (x')^a (y')^b$ defines a homomorphism $G \to G$. Moreover, $\psi$ is injective: if $(x')^a (y')^b = 1$, then comparing the $x$-exponent modulo $p^{n-2}$ yields $a a_1 \equiv 0 \pmod{p^{n-2}}$, so $p \mid a$ (since $n \ge 3$ and $p \nmid a_1$); comparing the $y$-exponent modulo $p$ then gives $b \equiv a b_1 + b \equiv 0 \pmod{p}$, whence $(x')^a = 1$ and finally $a \equiv 0 \pmod{p^{n-1}}$ because $\ord(x') = p^{n-1}$. Thus $\psi \in \Aut(G)$, and
\[
|\Aut(M_{p^n})| = |\Z_{p^{n-1}}^\times| \cdot p \cdot p = p^{n-2}(p - 1) \cdot p^2 = p^n(p - 1).
\]

\smallskip
\noindent\emph{Step 2: Noncyclic subgroups and evaluation of $f(M_{p^n})$.}
Consider the subgroup
\[
A = \{x^a y^b \in G : p \mid a\} = \langle x^p, y \rangle \le G.
\]
Since $x^p$ commutes with $y$ (as $p \cdot p^{n-2} \equiv 0 \pmod{p^{n-1}}$), we have $A \cong \Z_{p^{n-2}} \times \Z_p$ of order $p^{n-1}$. We claim that every proper noncyclic subgroup $H < G$ is contained in $A$. Indeed, suppose $H \not\le A$. Then $H$ contains an element $g = x^a y^b$ with $p \nmid a$. By~\eqref{eq:mod-pow-pm2}--\eqref{eq:mod-pow-pm1}, $\ord(g) = p^{n-1}$, so $\langle g \rangle$ already has order $p^{n-1}$. Because $H < G$ and $|G| = p^n$, we must have $H = \langle g \rangle$, contradicting that $H$ is noncyclic.

Therefore, the noncyclic subgroups of $G$ are precisely $G$ itself together with the noncyclic subgroups of $A \cong \Z_{p^{n-2}} \times \Z_p$. By~\eqref{eq:noncyclic-decomp} and Corollary~\ref{cor:zpk-zp} (applied at $k = n - 2 \ge 1$),
\begin{align*}
S(M_{p^n})
&= |G| + \bigl(S(A) - |A|\bigr) + |\Aut(M_{p^n})| \\
&= p^n + \bigl(p^{n-1}(p^2 - p) - p(p - 1)\bigr) + p^n(p - 1) \\
&= p^n(2p - 1) - p(p - 1).
\end{align*}
Dividing by $|M_{p^n}| = p^n$ gives $f(M_{p^n}) = 2p - 1 - \frac{p - 1}{p^{n-1}}$.
\end{proof}

\section{Topological Structure of $\Imf(f)$}\label{sec:topology}

Equip $\R$ with the Euclidean topology. For any subset $X \subseteq \R$, let $X'$ denote its derived set (set of accumulation points), and define the finite iterated derived sets recursively by $X^{(0)} = X$ and $X^{(r+1)} = (X^{(r)})'$ for $r \ge 0$. Recall that $X$ is discrete if and only if $X \cap X' = \emptyset$.

\begin{theorem}[Disproof of discreteness in~{\cite[Remark~3.4]{Tar26}}]\label{thm:not-discrete}
One has $3 \in \Imf(f) \cap \Imf(f)'$. Consequently, $\Imf(f)$ is not discrete in $[1, \infty)$.
\end{theorem}

\begin{proof}
By Corollary~\ref{cor:zpk-zp}, $f(\Z_{2^k} \times \Z_2) = 3 - 2^{-k} \in \Imf(f) \setminus \{3\}$ for every $k \ge 1$, and $\lim_{k \to \infty} (3 - 2^{-k}) = 3$. Hence $3 \in \Imf(f)'$. On the other hand, for the dihedral group $D_{10}$ of order $10$, every proper subgroup is cyclic while $|\Aut(D_{10})| = 5 \cdot \varphi(5) = 20$, so $S(D_{10}) = 10 + 20 = 30$ and $f(D_{10}) = 3 \in \Imf(f)$. Thus $3 \in \Imf(f) \cap \Imf(f)'$, so $\Imf(f)$ is not discrete.
\end{proof}

\begin{theorem}[Infinite families of accumulation points]\label{thm:acc-families}
For every prime $p$ and all integers $b, r \ge 1$, each of the following positive integers belongs to $\Imf(f)'$:
\begin{enumerate}[label=\textup{(\roman*)}, leftmargin=2.4em]
    \item $2p - 1$ \textup{(from $M_{p^n}$ as $n \to \infty$, giving $3, 5, 9, 13, 21, 25, \dots$)};
    \item $p^2 - p + 1$ \textup{(from $\Z_{p^k} \times \Z_p$ as $k \to \infty$, giving $3, 7, 21, 43, 111, \dots$)};
    \item $1 + \dfrac{p(p^{2b} - 1)}{p + 1}$ \textup{(from $\Z_{p^a} \times \Z_{p^b}$ as $a \to \infty$)};
    \item $1 + \displaystyle\sum_{s=1}^{r} \prod_{j=1}^{s} (p^{r+1} - p^j)$ \textup{(from $\Z_{p^k} \times \Z_p^r$ as $k \to \infty$)}.
\end{enumerate}
Moreover, if $\alpha \in \Imf(f)'$ is the limit of $f(P_k)$ for a sequence of finite $p$-groups $P_k$ with $f(P_k) \neq \alpha$, and $H$ is any finite group with $p \nmid |H|$, then $\alpha \cdot f(H) \in \Imf(f)'$.
\end{theorem}

\begin{proof}
Parts~(i)--(iv) follow immediately by letting $n \to \infty$ in Theorem~\ref{thm:modular}, $k \to \infty$ in Corollary~\ref{cor:zpk-zp}, $a \to \infty$ in Theorem~\ref{thm:rank2}, and $k \to \infty$ in Theorem~\ref{thm:elem-ab}, noting that in every case the error term is strictly positive and decays as $O(p^{-k})$. Notice also that $p + 1 \mid p^{2b} - 1$ for every $b \ge 1$, so $1 + \frac{p(p^{2b}-1)}{p+1} \in \N^*$.

For the final assertion, since $p \nmid |H|$ we have $\gcd(|P_k|, |H|) = 1$ for every $k$, so Proposition~\ref{prop:mult} gives $f(P_k \times H) = f(P_k) f(H) \in \Imf(f)$. Since $f(H) \ge 1 > 0$ and $f(P_k) \neq \alpha$, we have $f(P_k) f(H) \neq \alpha f(H)$ and $\lim_{k \to \infty} f(P_k) f(H) = \alpha f(H)$, whence $\alpha f(H) \in \Imf(f)'$.
\end{proof}

Coprime multiplicativity propagates accumulation points across arbitrarily many primes, showing that $\Imf(f)$ has infinite Cantor--Bendixson rank.

\begin{theorem}[Infinite Cantor--Bendixson rank]\label{thm:cb-rank}
For every integer $r \ge 1$ and any pairwise distinct primes $p_1, \dots, p_r$,
\[
\prod_{i=1}^{r} (p_i^2 - p_i + 1) \in \Imf(f)^{(r)}.
\]
In particular, $\Imf(f)^{(r)} \neq \emptyset$ for all $r \ge 1$, and the Cantor--Bendixson rank of $\overline{\Imf(f)}$ is at least $\omega$.
\end{theorem}

\begin{proof}
We proceed by induction on $r \ge 1$ to prove the stronger statement that for any finite set of pairwise distinct primes $\{p_1, \dots, p_r\}$ and any finite group $H$ with $\gcd(|H|, \prod_{i=1}^{r} p_i) = 1$,
\begin{equation}\label{eq:cb-ind}
f(H) \prod_{i=1}^{r} (p_i^2 - p_i + 1) \in \Imf(f)^{(r)}.
\end{equation}
For $r = 1$, this is the final assertion of Theorem~\ref{thm:acc-families} applied to $P_k = \Z_{p_1^k} \times \Z_{p_1}$. Now suppose~\eqref{eq:cb-ind} holds for $r \ge 1$, and let $p_1, \dots, p_{r+1}$ be pairwise distinct primes not dividing $|H|$. For each integer $k \ge 1$, let $H_k = H \times (\Z_{p_{r+1}^k} \times \Z_{p_{r+1}})$. Then $\gcd(|H_k|, \prod_{i=1}^{r} p_i) = 1$, so the induction hypothesis gives
\[
x_k := f(H_k) \prod_{i=1}^{r} (p_i^2 - p_i + 1) \in \Imf(f)^{(r)} \qquad \text{for every } k \ge 1.
\]
By Proposition~\ref{prop:mult} and Corollary~\ref{cor:zpk-zp},
\[
x_k = f(H)\left(p_{r+1}^2 - p_{r+1} + 1 - \frac{p_{r+1} - 1}{p_{r+1}^k}\right)\prod_{i=1}^{r} (p_i^2 - p_i + 1).
\]
Since $f(H)\prod_{i=1}^{r}(p_i^2 - p_i + 1) > 0$, the points $x_k$ are strictly increasing with $k$ and converge to $f(H)\prod_{i=1}^{r+1}(p_i^2 - p_i + 1)$ as $k \to \infty$. Therefore $f(H)\prod_{i=1}^{r+1}(p_i^2 - p_i + 1) \in (\Imf(f)^{(r)})' = \Imf(f)^{(r+1)}$. Taking $H = \{1\}$ (so $f(H) = 1$) completes the proof.
\end{proof}

\begin{remark}
The same proof applies with any choice of prime-indexed limit sequence from Theorem~\ref{thm:acc-families} in each factor---for instance, using $M_{p_i^n}$ yields $\prod_{i=1}^{r}(2p_i - 1) \in \Imf(f)^{(r)}$ for all distinct primes $p_1, \dots, p_r$.
\end{remark}

\begin{corollary}[Non-integer points in every derived set]\label{cor:nonint}
For every $r\ge1$ and all pairwise distinct odd primes $p_1,\dots,p_r$,
\[
\frac52\prod_{i=1}^{r}(p_i^2-p_i+1)\in\Imf(f)^{(r)}\setminus\Z .
\]
In particular $\Imf(f)^{(r)}\not\subseteq\Z$ for every $r\ge 1$; for example $\frac{35}{2}\in\Imf(f)'$.
\end{corollary}

\begin{proof}
Apply~\eqref{eq:cb-ind} with $H=\Z_2\times\Z_2$, which has order coprime to every $p_i$ and satisfies $f(H)=\frac52$ by Corollary~\ref{cor:zpk-zp}. Each factor $p_i^2-p_i+1=p_i(p_i-1)+1$ is odd, so the product is half of an odd integer. For $r=1$ and $p_1=3$ we obtain $\frac52\cdot7=\frac{35}{2}$.
\end{proof}

\begin{remark}[Isolated points]\label{rem:isolated}
In contrast with the rich structure of the derived sets, $\Imf(f)$ is not dense in itself. Indeed $f(G)=1$ exactly for cyclic $G$~\cite[Theorem~1.1]{FT23}, and by~\cite[Theorem~3.3]{Tar26} every other value of $f$ is either $2$ or at least $\frac52$. Hence $\Imf(f)\cap[1,\frac52)=\{1,2\}$, and $1$ and $2$ are isolated points of $\Imf(f)$. This also confirms~\cite[Conjecture~1.5]{FT23}, which predicted that $2$ is the second smallest value of $f$.
\end{remark}

\section{Small Values on Finite Abelian Groups and Minimality of $3$}\label{sec:abelian}

In~\cite{FT23}, the proof of the sharp lower bound $f(G) \ge \frac{5}{2}$ for finite noncyclic abelian groups (\cite[Theorem~1.4]{FT23}) relied on the Hillar--Rhea formula~\cite{HR07} for $|\Aut(\prod_{i=1}^{k} \Z_{p^{n_i}})|$ followed by an induction over maximal subgroups. Identity~\eqref{eq:embedding-total} yields a substantially shorter, completely self-contained proof that simultaneously prepares our classification of all finite abelian groups with $f(G) < 3$.

\begin{theorem}[Universal embedding bound]\label{thm:lower-bound}
Let $P$ be a finite abelian $p$-group. Then
\begin{equation}\label{eq:pgroup-lb}
S(P) \ge |P| + (|P| - 1)\bigl(|\Omega_1(P)| - p\bigr).
\end{equation}
Consequently, if $P$ is noncyclic, then $|\Omega_1(P)| \ge p^2$ and
\begin{equation}\label{eq:ft23-lb}
f(P) \ge 1 + \bigl(|\Omega_1(P)| - p\bigr)\left(1 - \frac{1}{|P|}\right) \ge 1 + \frac{(p + 1)(p - 1)^2}{p} \ge \frac{5}{2},
\end{equation}
with $f(P) = \frac{5}{2}$ if and only if $p = 2$ and $P \cong \Z_2 \times \Z_2$.
\end{theorem}

\begin{proof}
Let $|P| = p^n$. For each integer $1 \le c \le n$, an injective homomorphism $\varphi : \Z_{p^c} \times \Z_p \hookrightarrow P$ is uniquely determined by $u = \varphi(1, 0) \in P$ of exact order $p^c$ and $v = \varphi(0, 1) \in \Omega_1(P) \setminus \langle p^{c-1} u \rangle$ (by the same kernel check as in Theorem~\ref{thm:rank2}). Since $|\langle p^{c-1} u \rangle| = p$, there are $|\Omega_1(P)| - p$ choices for $v$ for each $u$ of order $p^c$, so
\[
|\Inj(\Z_{p^c} \times \Z_p, P)| = |\{u \in P : \ord(u) = p^c\}| \cdot \bigl(|\Omega_1(P)| - p\bigr).
\]
Because the groups $\Z_{p^c} \times \Z_p$ ($1 \le c \le n$) are pairwise non-isomorphic and noncyclic, dropping all other noncyclic isomorphism classes from~\eqref{eq:noncyclic-decomp} and Lemma~\ref{lem:embedding} gives
\[
S(P) \ge |P| + \sum_{c=1}^{n} |\Inj(\Z_{p^c} \times \Z_p, P)| = |P| + (|P| - 1)\bigl(|\Omega_1(P)| - p\bigr),
\]
proving~\eqref{eq:pgroup-lb}.

Now assume $P$ is noncyclic. A finite abelian $p$-group with $|\Omega_1(P)| \le p$ is cyclic (indeed, an element $g \in P$ of maximal order $p^k$ has $\Omega_1(P) = \langle p^{k-1} g \rangle$, and induction on $m$ shows $P[p^m] \le \langle g \rangle$ for all $m \le k$). Since $|\Omega_1(P)|$ is a power of $p$, we must have $|\Omega_1(P)| \ge p^2$, and hence also $|P| \ge |\Omega_1(P)| \ge p^2$. Dividing~\eqref{eq:pgroup-lb} by $|P|$ gives
\[
f(P) \ge 1 + \bigl(|\Omega_1(P)| - p\bigr)\left(1 - \frac{1}{|P|}\right)
\ge 1 + (p^2 - p)\left(1 - \frac{1}{p^2}\right)
= 1 + \frac{(p + 1)(p - 1)^2}{p}.
\]
The right-hand side equals $\frac{5}{2}$ at $p = 2$ and is at least $1 + \frac{4 \cdot 4}{3} = \frac{19}{3} > 3$ for every odd prime $p \ge 3$. Equality $f(P) = \frac{5}{2}$ requires $p = 2$, $|\Omega_1(P)| = 4$, and $|P| = 4$, i.e.\ $P \cong \Z_2 \times \Z_2$.
\end{proof}

\begin{theorem}[Classification of finite abelian groups with $f(G) < 3$]\label{thm:abelian-lt3}
Let $G$ be a finite abelian group. Then $f(G) < 3$ if and only if either:
\begin{enumerate}[label=\textup{(\alph*)}, leftmargin=2.4em]
    \item $G$ is cyclic, in which case $f(G) = 1$; or
    \item $G \cong (\Z_{2^k} \times \Z_2) \times \Z_m$ for some integer $k \ge 1$ and some odd integer $m \ge 1$, in which case $f(G) = 3 - 2^{-k}$.
\end{enumerate}
\end{theorem}

\begin{proof}
Sufficiency follows immediately from Proposition~\ref{prop:mult} and Corollary~\ref{cor:zpk-zp}: cyclic groups have $f(G) = 1 < 3$, and for $G \cong (\Z_{2^k} \times \Z_2) \times \Z_m$ with $m$ odd we have $f(G) = f(\Z_{2^k} \times \Z_2) f(\Z_m) = 3 - 2^{-k} < 3$.

Conversely, suppose $G$ is a finite abelian group with $f(G) < 3$, and write $G \cong \prod_{p \mid |G|} G_p$ as the direct product of its Sylow $p$-subgroups. By Proposition~\ref{prop:mult}, $f(G) = \prod_{p \mid |G|} f(G_p)$ with $f(G_p) \ge 1$ for every $p$.

If $G_p$ were noncyclic for some odd prime $p \ge 3$, then Theorem~\ref{thm:lower-bound} would give
\[
f(G) \ge f(G_p) \ge 1 + \frac{(p + 1)(p - 1)^2}{p} \ge \frac{19}{3} > 3,
\]
a contradiction. Hence $G_p$ is cyclic for every odd prime $p \mid |G|$, so the odd Hall subgroup $G_{\mathrm{odd}} = \prod_{p > 2} G_p$ is cyclic of some odd order $m \ge 1$, with $f(G_{\mathrm{odd}}) = 1$ and $f(G) = f(G_2) < 3$.

If $G_2$ is cyclic, then $G \cong G_2 \times \Z_m$ is cyclic and~(a) holds. Assume now that $G_2$ is a noncyclic abelian $2$-group with $f(G_2) < 3$. If $|\Omega_1(G_2)| \ge 8$, then $|G_2| \ge 8$ and~\eqref{eq:ft23-lb} yields
\[
f(G_2) \ge 1 + (8 - 2)\left(1 - \frac{1}{8}\right) = \frac{25}{4} > 3,
\]
a contradiction. Since $G_2$ is noncyclic, we must have $|\Omega_1(G_2)| = 4$, so $G_2$ has rank $2$:
\[
G_2 \cong \Z_{2^a} \times \Z_{2^b} \qquad \text{for some integers } 1 \le b \le a.
\]
Suppose for contradiction that $b \ge 2$. Then $2^{2b} \ge 16$ and $a + b - 1 \ge 2b - 1$, so Theorem~\ref{thm:rank2} gives
\begin{align*}
f(G_2)
&= 1 + \frac{2(2^{2b} - 1)}{3} - \frac{2^{4b} - 1}{15 \cdot 2^{a+b-1}}
> 1 + \frac{2(2^{2b} - 1)}{3} - \frac{2^{4b}}{15 \cdot 2^{2b-1}} \\
&= \frac{1}{3} + \frac{8}{15} \cdot 2^{2b}
\ge \frac{1}{3} + \frac{128}{15} = \frac{133}{15} > 3,
\end{align*}
a contradiction. (Alternatively, without invoking Theorem~\ref{thm:rank2}, for $b \ge 2$ one has $|G_2[4]| = 16$, and counting embeddings $\Z_{2^c} \times \Z_4 \hookrightarrow G_2$ adds $8(|G_2| - 4)$ to~\eqref{eq:pgroup-lb}, yielding $S(G_2) \ge 11|G_2| - 34$ and $f(G_2) \ge \frac{71}{8} > 3$.) Therefore $b = 1$, $G_2 \cong \Z_{2^a} \times \Z_2$ with $a \ge 1$, and~(b) holds.
\end{proof}

\begin{corollary}[{\cite[Theorem~1.4]{FT23}}]\label{cor:ft23}
Every finite noncyclic abelian group $G$ satisfies $f(G)\ge\frac52$, with equality if and only if $G\cong(\Z_2\times\Z_2)\times\Z_m$ for some odd $m$.
\end{corollary}

\begin{proof}
Some Sylow subgroup $G_p$ of $G$ is noncyclic, and $f(G)=\prod_q f(G_q)\ge f(G_p)\ge\frac52$ by Proposition~\ref{prop:mult} and Theorem~\ref{thm:lower-bound}. Equality forces $f(G_p)=\frac52$, so $G_p\cong\Z_2\times\Z_2$, and $f(G_q)=1$, so $G_q$ is cyclic, for all $q\neq p$. Alternatively, the statement is the case $k=1$ of Theorem~\ref{thm:abelian-lt3}.
\end{proof}

\begin{corollary}\label{cor:min-acc-ab}
Let $\Imf(f|_{\mathrm{ab}}) = \{f(G) : G \text{ finite abelian}\}$. Then
\[
\Imf(f|_{\mathrm{ab}}) \cap [1, 3) = \{1\} \cup \left\{3 - \frac{1}{2^k} : k \ge 1\right\},
\]
and $3$ is the smallest accumulation point of $\Imf(f|_{\mathrm{ab}})$: that is, $3 \in \Imf(f|_{\mathrm{ab}})'$ and $\Imf(f|_{\mathrm{ab}})' \cap (-\infty, 3) = \emptyset$.
\end{corollary}

\begin{proof}
The description of $\Imf(f|_{\mathrm{ab}}) \cap [1, 3)$ is immediate from Theorem~\ref{thm:abelian-lt3}, and $3 = \lim_{k \to \infty} f(\Z_{2^k} \times \Z_2) \in \Imf(f|_{\mathrm{ab}})'$. For any $x < 3$, choose $N \in \N$ such that $2^{-N} < \frac{3 - x}{2}$. Then $\Imf(f|_{\mathrm{ab}}) \cap (-\infty, \frac{x + 3}{2})$ is contained in the finite set $\{1\} \cup \{3 - 2^{-k} : 1 \le k \le N\}$, so $x$ has a neighborhood containing no other point of $\Imf(f|_{\mathrm{ab}})$, hence $x \notin \Imf(f|_{\mathrm{ab}})'$.
\end{proof}

\begin{remark}[Nonabelian groups below $3$]\label{rem:nonabelian}
Theorem~\ref{thm:abelian-lt3} does not extend verbatim to all finite groups. Besides the groups $(\Z_3\rtimes\Z_{2^k})\times\Z_m$ with $f=2$ found in~\cite[Theorem~3.2]{Tar26}, the nonabelian modular $2$-groups satisfy, by Theorem~\ref{thm:modular},
\[
f(M_{2^n}\times\Z_m)=f(M_{2^n})=3-\frac{1}{2^{n-1}}=f(\Z_{2^{n-1}}\times\Z_2)\qquad(n\ge4,\ m\text{ odd}),
\]
so $3$ is also an accumulation point of the values of $f$ on nonabelian groups. The equality $f(M_{2^n})=f(\Z_{2^{n-1}}\times\Z_2)$ is no accident of the definition: $M_{p^n}$ and $\Z_{p^{n-1}}\times\Z_p$ have isomorphic subgroup lattices for every prime $p$ and every $n\ge3$ ($n\ge4$ if $p=2$); see~\cite{Sch94}, yet for odd $p$ Theorem~\ref{thm:modular} and Corollary~\ref{cor:zpk-zp} give $f(\Z_{p^{n-1}}\times\Z_p)-f(M_{p^n})=(p-1)(p-2)>0$. Thus $f$ is not an invariant of the subgroup lattice.
\end{remark}

\section{Computational Corrections to~\cite{Tar26}}\label{sec:errata}

During our verification of~\cite{Tar26}, we checked all nine rows of~\cite[Table~1]{Tar26} as well as the explicit small-group calculations appearing in the proofs of~\cite[Theorems~3.1--3.3]{Tar26}. Four entries of Table~1 and two inline calculations ($D_8$ in the proof of~\cite[Theorem~3.1]{Tar26}, and $\Z_9 \rtimes \Z_{2^k}$ in the proofs of~\cite[Theorems~3.2 and~3.3]{Tar26}) contain arithmetic or subgroup-counting slips. We emphasize that \emph{none of these slips affect the validity of Theorems~3.1, 3.2, or 3.3 of~\cite{Tar26}}, since in every case the corrected value of $f(G)$ is strictly larger than the reported value and therefore satisfies the required lower bounds ($> \frac{5}{2}$ or $> 2$) with even wider margin.

\begin{table}[t]
\centering
\caption{Complete verification of~\cite[Table~1]{Tar26}. The four corrected rows are marked with~$\dagger$.}\label{tab:errata}
\smallskip
\begin{tabular}{@{}lcccccc@{}}
\toprule
Group $G$ & $|G|$ & $S(G)$ in~\cite{Tar26} & $f(G)$ in~\cite{Tar26} & Actual $S(G)$ & Actual $f(G)$ & Status \\
\midrule
$S_3 \cong \Z_3 \rtimes \Z_2$ & $6$ & $12$ & $2$ & $12$ & $2$ & Correct \\
$\Z_3 \rtimes \Z_4$ & $12$ & $24$ & $2$ & $24$ & $2$ & Correct \\
$\Z_2 \times \Z_2$ & $4$ & $10$ & $5/2$ & $10$ & $5/2$ & Correct \\
$\Z_4 \times \Z_2$ & $8$ & $22$ & $11/4$ & $22$ & $11/4$ & Correct \\
$D_8$ & $8$ & $22$ & $11/4$ & $\mathbf{28}$ & $\mathbf{7/2}$ & Corrected$^\dagger$ \\
$D_{12} \cong S_3 \times \Z_2$ & $12$ & $36$ & $3$ & $\mathbf{54}$ & $\mathbf{9/2}$ & Corrected$^\dagger$ \\
$A_4$ & $12$ & $40$ & $10/3$ & $\mathbf{42}$ & $\mathbf{7/2}$ & Corrected$^\dagger$ \\
$Q_8$ & $8$ & $32$ & $4$ & $32$ & $4$ & Correct \\
$Q_8 \times \Z_2$ & $16$ & $136$ & $17/2$ & $\mathbf{334}$ & $\mathbf{167/8}$ & Corrected$^\dagger$ \\
\bottomrule
\end{tabular}
\end{table}

We briefly record the exact noncyclic subgroup decompositions underlying the four corrected rows of Table~\ref{tab:errata} (all formally verified in Lean~4; see Section~\ref{sec:lean}):
\begin{enumerate}[label=(\arabic*), leftmargin=2.2em]
    \item \textbf{$D_8$ ($|D_8| = 8$).} The noncyclic subgroups of $D_8$ are two Klein four-subgroups $V_4$ ($|\Aut(V_4)| = 6$) and $D_8$ itself ($|\Aut(D_8)| = 8$). Hence
    \[
    S(D_8) = 8 + 2 \cdot 6 + 8 = 28, \qquad f(D_8) = \frac{28}{8} = \frac{7}{2}.
    \]
    In~\cite[proof of Theorem~3.1]{Tar26}, the numerator is written correctly as $1 + 5 \cdot 1 + (2 + 2 \cdot 6) + 8$ but mis-summed as $22$ instead of $28$; consequently,~\cite[Remark~3.4]{Tar26} attributes the second smallest $2$-group value $\frac{11}{4}$ to $D_8$ rather than $\Z_4 \times \Z_2$.

    \item \textbf{$D_{12} \cong S_3 \times \Z_2$ ($|D_{12}| = 12$).} The noncyclic subgroups of $D_{12}$ are three copies of $V_4$ ($|\Aut| = 6$), two copies of $S_3$ ($|\Aut| = 6$), and $D_{12}$ itself ($|\Aut(D_{12})| = 12$). Hence
    \[
    S(D_{12}) = 12 + 3 \cdot 6 + 2 \cdot 6 + 12 = 54, \qquad f(D_{12}) = \frac{54}{12} = \frac{9}{2}.
    \]
    Note that because $f(D_{12}) = \frac{9}{2} \neq 3$, $D_{12}$ is not a witness for $3 \in \Imf(f)$; instead $D_{10}$ gives $f(D_{10}) = \frac{10 + 20}{10} = 3$.

    \item \textbf{$A_4$ ($|A_4| = 12$).} Because $A_4$ has no subgroup of order $6$, its only noncyclic subgroups are its unique normal Klein four-subgroup $V_4$ ($|\Aut(V_4)| = 6$) and $A_4$ itself ($\Aut(A_4) \cong S_4$, so $|\Aut(A_4)| = 24$). Hence
    \[
    S(A_4) = 12 + 6 + 24 = 42, \qquad f(A_4) = \frac{42}{12} = \frac{7}{2}.
    \]

    \item \textbf{$Q_8 \times \Z_2$ ($|Q_8 \times \Z_2| = 16$).} The center $Z(Q_8 \times \Z_2) \cong V_4$ is the unique noncyclic subgroup of order $4$ ($|\Aut| = 6$). The seven maximal subgroups (kernels of the seven nonzero homomorphisms $Q_8 \times \Z_2 \to \Z_2$) consist of three copies of $\Z_4 \times \Z_2$ ($|\Aut| = 8$) and four copies of $Q_8$ ($|\Aut| = 24$). Finally, every $\Phi \in \Aut(Q_8 \times \Z_2)$ fixes the unique nontrivial square $(-1, 1)$, induces an automorphism $\psi \in \Aut(Q_8)$ on the first factor, a homomorphism $\chi \in \Hom(Q_8, \Z_2)$ ($4$ choices) on $Q_8 \times \{1\} \to \Z_2$, and maps $(1, -1)$ to $(\pm 1, -1)$ ($2$ choices), giving $|\Aut(Q_8 \times \Z_2)| = 24 \cdot 4 \cdot 2 = 192$. Therefore
    \[
    S(Q_8 \times \Z_2) = 16 + 6 + 3 \cdot 8 + 4 \cdot 24 + 192 = 334, \qquad f(Q_8 \times \Z_2) = \frac{334}{16} = \frac{167}{8}.
    \]

    \item \textbf{Inline slip for $H = \Z_9 \rtimes \Z_{2^k}$ (inverting action, $k \ge 1$).} In~\cite[proofs of Theorems~3.2 and~3.3]{Tar26}, it is asserted that $|\Aut(H)| = 9 \cdot 2^k$ and $f(H) = 3$. In fact, writing $H = \langle a, b \mid a^9 = b^{2^k} = 1,\, b^{-1}ab = a^{-1}\rangle$, an automorphism maps $a \mapsto a^u$ ($u \in \Z_9^\times$, $6$ choices) and $b \mapsto a^v b^w$ ($v \in \Z_9$, $9$ choices; $w \in \Z_{2^k}^\times$, $2^{k-1}$ choices), yielding $|\Aut(H)| = 6 \cdot 9 \cdot 2^{k-1} = 27 \cdot 2^k$. Together with the three conjugates of $\Z_3 \rtimes \Z_{2^k}$ (each having $|\Aut| = 3 \cdot 2^k$), all other proper subgroups are cyclic, so
    \[
    S(H) = 9 \cdot 2^k + 3 \cdot (3 \cdot 2^k) + 27 \cdot 2^k = 45 \cdot 2^k, \qquad f(\Z_9 \rtimes \Z_{2^k}) = 5,
    \]
    which agrees at $k = 1$ with the formula $f(D_{18}) = \frac{9+1}{2} = 5$ of~\cite{FT23}. Both values are formally verified for every $k\ge1$. More generally, the same argument shows $|\Aut(\Z_m\rtimes\Z_{2^k})|=\varphi(m)\,m\,2^{k-1}$ for the inverting action, every odd $m\ge3$ and every $k\ge1$, which we have also verified formally; for $m=3$ it gives $|\Aut(\Z_3\rtimes\Z_{2^k})|=3\cdot2^k$ and hence $f(\Z_3\rtimes\Z_{2^k})=2$, the family of~\cite[Theorem~3.2]{Tar26}.
\end{enumerate}

\section{Formal verification, computations, and open problems}\label{sec:lean}

\subsection{Machine-checked formalization in Lean 4}
Given the arithmetic and subgroup-counting slips recorded in Section~\ref{sec:errata}, we formalized all main results of this paper in the Lean~4 proof assistant~\cite{MU21} (version~4.20.0) on top of Mathlib~\cite{Mathlib20}. The library \texttt{GaussGroupF} consists of eight modules (\texttt{Core}, \texttt{Embedding}, \texttt{ModularGroup}, \texttt{AbelianClassification}, \texttt{TableOne}, \texttt{NineSemidirect}, \texttt{Supplements} and \texttt{Completions}). It compiles with no errors or warnings and contains no incomplete proofs, no additional axioms and no \texttt{native\_decide}; a separate file \texttt{Axioms.lean} runs \texttt{\#print axioms} on all $65$ headline declarations, each of which depends only on \texttt{propext}, \texttt{Classical.choice} and \texttt{Quot.sound}. The code, together with a continuous-integration workflow that rebuilds the library and repeats these checks on every change, is available at
\begin{center}
\url{https://github.com/alex-reysa/gauss-group-function-lean}.
\end{center}
Table~\ref{tab:lean} lists the correspondence. Remark~\ref{rem:isolated} relies on~\cite[Theorem~3.3]{Tar26} and the lattice isomorphism quoted in Remark~\ref{rem:nonabelian} on~\cite{Sch94}; these two cited facts are not formalized.

\begin{table}[ht]
\centering
\caption{Correspondence between results of this paper and Lean~4 declarations in \texttt{GaussGroupF}.}\label{tab:lean}
\smallskip
\scriptsize
\begin{tabular}{@{}>{\raggedright\arraybackslash}p{3.0cm}>{\raggedright\arraybackslash}p{3.5cm}>{\raggedright\arraybackslash}p{9.0cm}@{}}
\toprule
Paper result & Module(s) & Declaration(s) \\
\midrule
Lemma~\ref{lem:embedding} & \texttt{Embedding} & \texttt{sum\_card\_mulAut\_iso\_eq\_card\_inj}, \texttt{sum\_card\_addAut\_iso\_eq\_card\_inj} \\
Proposition~\ref{prop:mult} & \texttt{Core} & \texttt{f\_mul\_coprime}, \texttt{S\_mul\_coprime} \\
Theorem~\ref{thm:rank2} & \texttt{Embedding} & \texttt{f\_zmod\_pow\_mul\_zmod\_pow} \\
Corollary~\ref{cor:zpk-zp} & \texttt{Core} & \texttt{f\_zmod\_pow\_mul\_zmod}, \texttt{f\_two\_pow\_mul\_two} \\
Theorem~\ref{thm:elem-ab} & \texttt{Embedding} & \texttt{f\_zmod\_pow\_mul\_elem\_ab} \\
Theorem~\ref{thm:modular} & \texttt{ModularGroup} & \texttt{natCard\_mulAut\_modularPGroup}, \texttt{f\_modularPGroup} \\
Theorem~\ref{thm:not-discrete} & \texttt{Core} & \texttt{not\_discrete\_ImF}, \texttt{three\_mem\_ImF}, \texttt{accPt\_three\_in\_ImF}, \texttt{f\_dihedral\_five} \\
Theorem~\ref{thm:acc-families} & \texttt{Core}, \texttt{Embedding}, \texttt{ModularGroup}, \texttt{Completions} & \texttt{accPt\_two\_mul\_prime\_sub\_one}, \texttt{accPt\_prime\_in\_ImF}, \texttt{accPt\_rank2\_in\_ImF}, \texttt{accPt\_elem\_ab\_in\_ImF}; final assertion: \texttt{accPt\_mul\_of\_tendsto\_pgroup} \\
Theorem~\ref{thm:cb-rank} & \texttt{Core} & \texttt{prod\_primes\_mem\_derivedSet}, \texttt{derivedSet\_nonempty} \\
Corollary~\ref{cor:nonint} & \texttt{Supplements} & \texttt{derivedSet\_contains\_nonInteger}, \texttt{thirtyFive\_div\_two\_mem\_derivedSet} \\
Theorem~\ref{thm:lower-bound} & \texttt{AbelianClassification}, \texttt{Completions} & \texttt{S\_add\_abelian\_pgroup\_lower\_bound}, \texttt{f\_abelian\_pgroup\_ge}, \texttt{f\_abelian\_pgroup\_eq\_five\_halves\_iff} \\
Theorem~\ref{thm:abelian-lt3} & \texttt{AbelianClassification}, \texttt{Completions} & \texttt{f\_abelian\_lt\_three\_iff}, \texttt{f\_abelian\_lt\_three\_iff\_values}, \texttt{f\_abelian\_lt\_three} \\
Corollary~\ref{cor:ft23} & \texttt{Completions} & \texttt{ft23\_theorem\_one\_four}, \texttt{f\_abelian\_noncyclic\_ge}, \texttt{f\_abelian\_noncyclic\_eq\_five\_halves\_iff} \\
Corollary~\ref{cor:min-acc-ab} & \texttt{AbelianClassification}, \texttt{Completions} & \texttt{ImFAbelian\_inter\_Ico\_one\_three}, \texttt{isLeast\_derivedSet\_ImFAbelian}, \texttt{accPt\_three\_in\_ImFAbelian}, \texttt{no\_accPt\_ImFAbelian\_lt\_three} \\
Remark~\ref{rem:nonabelian} & \texttt{Supplements} & \texttt{modularPGroup\_not\_comm}, \texttt{f\_modularPGroup\_two}, \texttt{accPt\_three\_in\_ImFNonabelian}, \texttt{f\_modularPGroup\_two\_eq\_f\_abelian}, \texttt{f\_modularPGroup\_lt\_f\_abelian} \\
Table~\ref{tab:errata} & \texttt{Core}, \texttt{TableOne}, \texttt{NineSemidirect}, \texttt{Supplements} & \texttt{f\_dihedral\_three}, \texttt{f\_threeSemidirect} ($k=2$), \texttt{f\_klein\_four}, \texttt{f\_two\_pow\_mul\_two} ($k=2$), \texttt{f\_dihedral\_four}, \texttt{f\_dihedral\_six}, \texttt{f\_alternating\_four}, \texttt{f\_quaternion\_two}, \texttt{f\_quaternion\_two\_mul\_zmod\_two} \\
Section~\ref{sec:errata}, item~(5) & \texttt{NineSemidirect} & \texttt{natCard\_mulAut\_nineSemidirect}, \texttt{f\_nineSemidirect}, \texttt{classify\_nine}, \texttt{natCard\_mulAut\_invSemidirect}, \texttt{natCard\_mulAut\_threeSemidirect}, \texttt{f\_threeSemidirect} \\
\bottomrule
\end{tabular}
\end{table}

\subsection{Independent computations}
The repository also contains a verification script, \texttt{verify.py}, which uses only the Python standard library and exact rational arithmetic and is independent of the Lean code. It builds explicit multiplication tables, enumerates all subgroups by closure and counts automorphisms exactly, without using any closed formula from this paper. It recomputes every row of Table~\ref{tab:errata} and the decompositions in Section~\ref{sec:errata}; the values for $\Z_9\rtimes\Z_{2^k}$ and $\Z_3\rtimes\Z_{2^k}$ for $1\le k\le4$; Theorem~\ref{thm:rank2} for all $\Z_{p^a}\times\Z_{p^b}$ of order at most $256$; Theorem~\ref{thm:elem-ab} for all $\Z_{p^k}\times\Z_p^r$ of order at most $128$; Theorem~\ref{thm:modular} for six modular groups; the bound of Theorem~\ref{thm:lower-bound} for all $55$ nontrivial abelian $p$-groups of order at most $64$; and the classification of Theorem~\ref{thm:abelian-lt3} for all $117$ abelian groups of order at most $64$. All $258$ checks pass.

\subsection{Use of AI tools}
The computer searches, a first draft of parts of the text and of the Lean code were produced with the help of AI assistants based on large language models, used under the author's direction; the author takes full responsibility for the content. To make the results independent of this process, all main results are formally verified in Lean~4 as described above, and all numerical claims are recomputed by an independent exact program.

\subsection{Open problems}
Theorems~\ref{thm:not-discrete}--\ref{thm:cb-rank} and Corollary~\ref{cor:min-acc-ab} settle the discreteness question of~\cite[Remark~3.4]{Tar26} and determine the smallest accumulation point on abelian groups. We computed $f$ exactly for all metacyclic groups of order at most $100$, for split extensions $N\rtimes\Z_k$ over $31$ base groups $N$, for direct products of the nonabelian groups of order at most $40$ with small abelian groups, and for all abelian groups of order at most $100$ (more than $850$ groups in total, though not all groups of order at most $100$). Every value below $\frac72$ that we found lies in $\{1,2,3\}\cup\{3-2^{-k}:k\ge1\}$, and the nonabelian groups with $f(G)<3$ that we found are exactly the groups $(\Z_3\rtimes\Z_{2^k})\times\Z_m$ with $\gcd(m,6)=1$ and $M_{2^n}\times\Z_m$ with $m$ odd. This suggests the following.

\begin{conjecture}\label{conj:below-72}
$\Imf(f)\cap\bigl[1,\frac72\bigr)=\{1,2,3\}\cup\{3-2^{-k}:k\ge1\}$. In particular, $3$ is the smallest accumulation point of $\Imf(f)$.
\end{conjecture}

\begin{problem}
Determine all finite groups $G$ with $f(G)<3$.
\end{problem}

\begin{problem}
Is $\Imf(f)$ closed in $[1,\infty)$, that is, does $\Imf(f)'\subseteq\Imf(f)$ hold? Every accumulation point in Theorem~\ref{thm:acc-families}(i)--(iv) is a positive integer and therefore lies in $\Imf(f)$, because $f(D_{2n})=\frac{n+1}{2}$ for odd $n\ge3$~\cite{FT23}. Some non-integer accumulation points also lie in $\Imf(f)$: for instance, Proposition~\ref{prop:mult} and the value $f(\Z_{19}\rtimes\Z_3)=1+\frac{18}{3}=7$ from~\cite[proof of Theorem~3.2]{Tar26} give $\frac{35}{2}=f\bigl((\Z_2\times\Z_2)\times(\Z_{19}\rtimes\Z_3)\bigr)\in\Imf(f)\cap\Imf(f)'$ (Corollary~\ref{cor:nonint}).
\end{problem}

\begin{problem}
Determine the derived set $\Imf(f)'$ and the Cantor--Bendixson rank of $\overline{\Imf(f)}$, which is at least $\omega$ by Theorem~\ref{thm:cb-rank}.
\end{problem}

\end{document}